\documentclass{article}

\usepackage{amsthm,amssymb,booktabs,array}

\def\Res{\mathop{\operatorfont res}\nolimits}
\def\Pic{\mathop{\operatorfont Pic}\nolimits}
\def\supp{\mathop{\rm supp}\nolimits}

\usepackage{amsmath}

\newtheorem{theorem}{Theorem}[section]
\newtheorem{corollary}[theorem]{Corollary}
\newtheorem{proposition}[theorem]{Proposition}
\newtheorem{lemma}[theorem]{Lemma}

\theoremstyle{definition}

\newcommand{\ff}{{\mathbb F}}

\def\ord{\mathop{\operatorfont ord}\nolimits}
\def\wt{\mathop{\operatorfont wt}\nolimits}

\def\Res{\mathop{\operatorfont Res}\nolimits}

\usepackage{url}

\begin{document}


\title{An extension of the order bound for AG codes}

\author{Iwan Duursma and Radoslav Kirov}
%


\maketitle

\begin{abstract}
The most successful method to obtain lower bounds for the minimum distance of an algebraic geometric code is
the order bound, which generalizes the Feng-Rao bound. We provide a significant extension 
of the bound that improves the order bounds by Beelen and by Duursma and Park. We include
an exhaustive numerical comparison of the different bounds for $10168$ two-point codes 
on the Suzuki curve of genus $g=124$ over the field of $32$ elements. \\[1ex]
\end{abstract}

\section{Introduction}

To obtain lower bounds
for the minimum distance of algebraic geometric codes we follow \cite{DuuPar08} and exploit 
a relation between the minimum distance of algebraic geometric codes and the 
representation of divisor classes by differences of base point free divisors. \\

Let $X/\ff$ be an algebraic curve (absolutely irreducible, smooth,
projective) of genus $g$ over a finite field $\ff$. Let $\ff(X)$
be the function field of $X/\ff$. A nonzero rational function $f \in \ff(X)$
has divisor $(f) = \sum_{P \in X} \ord_P(f) P = (f)_0 - (f)_\infty,$ where the 
positive part $(f)_0$ gives the zeros of $f$ and their multiplicities, and the
negative part $(f)_\infty$ gives the poles of $f$ and their multiplicities. 
A divisor $D = \sum_P m_P P$ is principal if it is of the form $D=(f)$ for some
$f \in \ff(X).$ Two divisors $D$ and $D'$ are linearly equivalent if $D' = D + (f)$
for some $f \in \ff(X).$ Given a divisor $D$ on $X$ defined over $\ff$, let $L(D)$ 
denote the vector space over $\ff$ of nonzero functions $f \in \ff(X)$ for which
$(f) + D \geq 0$ together with the zero function. A point $P$ is a base point for
the linear system of divisors $\{ (f) + D : f \in L(D) \}$ if $(f) + D \geq P$ for
all $f \in L(D)$, that is to say if $L(D) = L(D-P).$ \\

We give the definition of an algebraic geometric code. For $n$ distinct rational points 
$P_1, \ldots, P_n$ on $X$ and for disjoint divisors $D=P_1+\cdots+P_n$ and $G$, the 
geometric Goppa code $C_L(D,G)$ is defined as the image of the map
\[
\alpha _L ~:~ L(G)~~\longrightarrow~~\ff^{\,n}, ~~f \mapsto (\,f(P_1), \ldots, f(P_n) \,).
\]
With the Residue theorem, the dual code $C_L(D,G)^\perp$ can be expressed in terms of 
differentials. Let $\Omega(X)$ be the module of rational differentials for $X/\ff$. For a given divisor
$E$ on $X$ defined over $\ff$, let 
$\Omega(E)$ denote the vector space over $\ff$ of nonzero differentials $\omega \in \Omega(X)$ for which
$(\omega) \geq E$ together with the zero differential. Let $K$ represent
the canonical divisor class. The geometric Goppa code $C_\Omega(D,G)$
is defined as the image of the map
\[
\alpha _\Omega ~:~ \Omega(G-D)~~\longrightarrow~~\ff^{\,n}, 
~~\omega \mapsto (\,\Res_{P_1}(\omega), \ldots, \Res_{P_n}(\omega) \,).
\]
The code $C_\Omega(D,G)$ is the dual code for the code $C_L(D,G).$ 
We use the following characterization of the minimum distance.

\begin{proposition} \cite[Proposition 2.1]{DuuPar08} \label{P:agd}
For the code $C_L(D,G)$, and for $C=D-G$, 
\[
d(C_L(D,G)) = \min \{ \deg A : 0 \leq A \leq D \;|\; L(A-C) \neq L(-C) \}. 
\]
For the code $C_\Omega(D,G)$, and for $C=G-K$,
\[
d(C_\Omega(D,G)) = \min \{ \deg A : 0 \leq A \leq D \;|\; L(A-C) \neq L(-C) \}.
\]
\end{proposition}
\begin{proof}
There exists a nonzero word in $C_L(D,G)$ with support in $A$, for $0 \leq A \leq D$,
if and only if $L(G-D+A)/L(G-D) \neq 0.$ There exists a nonzero word in $C_\Omega(D,G)$ with support in $A$,
for $0 \leq A \leq D$, if and only if $\Omega(G-A)/\Omega(G) \neq 0$ if and only if
$L(K-G+A) / L(K-G) \neq 0.$
\end{proof}

It is clear that in each case $d \geq \deg C$. The lower bound $d_{GOP} = \deg C$ is the Goppa designed
minimum distance of a code.
 
\section{Coset bounds}

For a point $P$ disjoint from $D$, consider the subcodes $C_L(D,G-P) \subseteq C_L(D,G)$ and
$C_\Omega(D,G+P) \subseteq C_\Omega(D,G).$ 

\begin{proposition} \cite[Proposition 3.5]{DuuPar08} \label{P:cb}
Let $A = \supp(c)$ be the support of a codeword $c = (c_P : P \in D)$, with $0 \leq A \leq D.$
For $C=D-G,$
\[
c \in C_L(D,G) \backslash C_L(D,G-P) ~\Rightarrow~ L(A-C) \neq L(A-C-P).
\]
For $C=G-K,$
\[
c \in C_\Omega(D,G) \backslash C_\Omega(D,G+P) ~\Rightarrow~ L(A-C) \neq L(A-C-P).
\]
\end{proposition}
\begin{proof}
There exists a word in $C_L(D,G) \backslash C_L(D,G-P)$ with support in $A$, for $0 \leq A \leq D$,
if and only if $C_L(D,G-(D-A)) \neq C_L(D,G-P-(D-A))$ only if $L(G-D+A) \neq L(G-D+A-P).$ There exists a word in 
$C_\Omega(D,G) \backslash C_\Omega(D,G+P)$ with support in $A$,
for $0 \leq A \leq D$, if and only if $C_\Omega(A,G) \neq C_\Omega(A,G+P)$ only if $\Omega(G-A) \neq \Omega(G-A+P)$,
which can be expressed as $L(K-G+A) \neq L(K-G+A-P).$
\end{proof}

The order bound for the minimum distance of an algebraic geometric code is motivated by the decoding procedures
in \cite{FenRao93}, \cite{Duu93} and combines estimates for the weight of a word 
$c \in C_L(D,G) \backslash$ $C_L(D,G-P)$, or $c \in C_\Omega(D,G) \backslash$ $C_\Omega(D,G+P)$.
The basic version (often referred to as the simple or first order bound \cite{HoeLinPel98}, \cite{CamFarMun00}) 
takes the form
\begin{align*}
&d(C_L(D,G)) = \min_{i \geq 0} ( \, \min \{ \wt(c) : c \in C_L(D,G-iP) \backslash C_L(D,G-(i+1)P) \} \, ), \\
&d(C_\Omega(D,G)) = \min_{i \geq 0} ( \, \min \{ \wt(c) : c \in C_\Omega(D,G+iP) \backslash C_L(D,G+(i+1)P) \} \, ).
\end{align*}
The order bound makes it possible to use separate estimates for different subsets of codewords. The bound is
successful if for each subset, i.e. for each $i$, we can find an estimate that is better than a uniform lower 
bound for all codewords. Methods that provide uniform lower bounds include the Goppa designed minimum distance 
and bounds of floor type \cite{MahMat06}, \cite{LunMcc06}, \cite{GunStiTas09}, \cite{DuuPar08}. It follows 
from the Singleton bound that the minimum distance of an algebraic geometric 
code can not exceed its designed minimum distance by more than $g$, where $g$ is the genus of the curve.
This implies that the minimum in the order bound occurs for $i \in \{0,\ldots,g\}.$ \\

For a curve $X$ defined over the field $\ff$, let $\Pic(X)$ be the group of divisor classes. 
Let $\Gamma = \{ A : L(A) \neq 0 \}$ be the semigroup of effective divisor classes. For a divisor
class $C$, define 
\begin{align*}
\Gamma(C) &= \{ A : L(A) \neq 0 \,\text{ and }\, L(A-C) \neq 0 \}, \\
\Gamma^\ast(C) &= \{ A : L(A) \neq 0 \,\text{ and }\, L(A-C) \neq L(-C) \}.
\end{align*} 
The semigroup $\Gamma(C)$ has the property that $A+E \in \Gamma(C)$ whenever $A \in \Gamma(C)$ and $E \in \Gamma.$ 
With the extra structure $\Gamma(C)$ is a semigroup ideal. Similar for $\Gamma^\ast(C).$ \\

For a suitable choice of divisor class $C$, the subsets of coordinates $A$ that support a codeword
in an algebraic geometric code $C_L(D,G)$ or $C_\Omega(D,G)$ belong to the semigroup ideals
$\Gamma^\ast(C) \subseteq \Gamma(C).$
\begin{align*}
c \in C_L(D,G) \backslash \{0\} &~\Rightarrow~ A = \supp(c) \in \Gamma^\ast(C) \subseteq \Gamma(C), \\
c \in C_\Omega(D,G) \backslash \{0\} &~\Rightarrow~ A = \supp(c) \in \Gamma^\ast(C) \subseteq \Gamma(C).
\end{align*}
Following \cite{DuuPar08}, our approach from here on will be to estimate the minimal degree of a divisor $A \in \Gamma^\ast(C)$.
For that purpose we no longer need to refer to the codes $C_L(D,G)$ or $C_\Omega(D,G)$ after we
choose $C=D-G$ or $C=G-K$, respectively.
We write ${\mathcal C}(C)$ to refer to any code with
designed minimum support $D-G$ or $G-K$ in the divisor class $C$. In this short paper, we restrict ourselves to the
case $\deg C > 0$ (i.e. to codes with positive Goppa designed minimum distance), 
so that $L(-C) = 0$, and $\Gamma^\ast(C) = \Gamma(C).$ The case $\deg C \leq 0$ 
is handled with a straightforward modification similar to that used in \cite{DuuPar08}. \\

To apply the order bound argument we restate Proposition \ref{P:cb} in terms of semigroup ideals.
For a given point $P \in X$, let $\Gamma_P = \{ A : L(A) \neq L(A-P) \}$ be the semigroup of effective 
divisor classes with no base point at $P$. For a divisor class $C$ and for a point $P$, define
the semigroup ideal
\[
\Gamma_P(C) = \{ A : L(A) \neq L(A-P) \,\text{ and }\, L(A-C) \neq L(A-C-P) \}.
\]
The implications in Proposition \ref{P:cb} become
\begin{align*}
c \in C_L(D,G) \backslash C_L(D,G-P) &~\Rightarrow~ A = \supp(c) \in \Gamma_{P}(C), \\
c \in C_\Omega(D,G) \backslash C_\Omega(D,G+P) &~\Rightarrow~ A = \supp(c) \in \Gamma_{P}(C). 
\end{align*}
Let $\Delta_P(C)$ be the complement of $\Gamma_P(C)$ in $\Gamma_P$,
\[
\Delta_P(C) = \{ A : L(A) \neq L(A-P) \,\text{ and }\, L(A-C) = L(A-C-P) \}.
\]

\begin{theorem} (Duursma-Park \cite{DuuPar08}) \label{T:cbdiv}
Let $\{ A_1 \leq A_2 \leq \cdots \leq A_w \} \subset \Delta_P(C)$ be
a sequence of divisors with $A_{i+1} \geq A_i + P$, for $i=1,\ldots,w-1.$
Then $\deg A \geq w$, for every divisor $A \in \Gamma_P(C)$ with support disjoint
from $A_w-A_1.$
\end{theorem}
\begin{proof}
A more general version is proved in the next section.
\end{proof}

\section{Order bounds}

Following \cite{Bee07FF}, we use the order bound with a sequence of points $\{ Q_i : i \geq 0 \}.$
For $R_i = Q_0 + Q_1 + \cdots + Q_{i-1},$ $i \geq 0,$
\begin{align*}
c \in C_L(D,G-R_i) \backslash C_L(D,G-R_i-Q_i) &~\Rightarrow~ A = \supp(c) \in \Gamma_{Q_i}(C+R_i), \\
c \in C_\Omega(D,G+R_i) \backslash C_\Omega(D,G+R_i+Q_i) &~\Rightarrow~ A = \supp(c) \in \Gamma_{Q_i}(C+R_i). 
\end{align*} 
Let ${\mathcal C}(C,S)$ refer to any code with designed minimum support in the divisor class $C$ and 
with divisor $D$ disjoint from $S$. Repeated application of Theorem \ref{T:cbdiv} gives the following bound for
the minimum distance. For comparison we state the Beelen bound \cite[Theorem 7, Remark 5]{Bee07FF} in a similar form.

\begin{corollary} (Duursma-Park bound $d_{DP}$ for the minimum distance)
For a code ${\mathcal C}(C,S)$, let $\{ Q_i : i \geq 0 \}$ be a sequence of points in $S$, such that for each 
$i \geq 0$ there exists a sequence $\{ A_1 \leq A_2 \leq \cdots \leq A_w \} \subset \Delta_{Q_i}(C+R_i)$,
for $R_i = Q_0 + Q_1 + \cdots + Q_{i-1},$ such that $A_{j+1} \geq A_j + Q_i$, for $j=1,\ldots,w-1,$ and
$A_w-A_1$ has support in $S$. Then ${\mathcal C}(C,S)$ has minimum distance at least $w$.
\end{corollary}

\begin{corollary} (Beelen bound $d_B$ for the minimum distance)
For a code ${\mathcal C}(C,S)$, let $\{ Q_i : i \geq 0 \}$ be a sequence of points in $S$, such that for each 
$i \geq 0$ there exists a sequence $\{ A_1 \leq A_2 \leq \cdots \leq A_w \} \subset \Delta_{Q_i}(C+R_i)$,
for $R_i = Q_0 + Q_1 + \cdots + Q_{i-1},$ such that $A_{j+1} - A_j \in \{ k Q_i : k > 0 \}$, for $j=1,\ldots,w-1.$ 
Then ${\mathcal C}(C,S)$ has minimum distance at least $w$.
\end{corollary}

The Beelen bound $d_B$ has a weaker version (which we will denote by $d_{B_0}$) that assumes that moreover 
$A_1 \in \{ k Q_i : k \geq 0 \}$.
The simple or first order bound uses the further specialization that $Q_i = P$ for $i \geq 0.$
For $G = K+C,$ such that $\deg C > 0$, and for $D$ disjoint from $P$, the simple order bound
becomes
\[
d(C_L(D,G)^\perp) \geq \min_{i \geq 0} \# (\Delta_P(C+iP) \cap \{ jP : j \geq 0 \}).
\]
For $G=mP$ this is the original Feng-Rao bound. \\

The purpose of the order bound is to improve on uniform bounds such as the floor bound.
In rare occasions the Beelen bound is less than bounds of floor type \cite{Bee07FF}, \cite{DuuPar08} (this is the case
for a single code in Table \ref{Table:comparison}). Compared to the Beelen bound, the ABZ order bound $d_{ABZ'}$ 
\cite[Theorem 6.6]{DuuPar08} allows $A_{j+1} - A_j \not \in \{ k Q_i : k > 0 \}$ for a single $j$ in the range $j=1,\ldots,w-1.$ With
that modification the ABZ order bound always is at least the ABZ floor bound $d_{ABZ}$ \cite[Theorem 2.4]{DuuPar08} which
is the best known bound of floor type. In general, $d_{B_0} \leq d_B \leq d_{ABZ'} \leq d_{DP}$ and
$d_{ABZ} \leq d_{ABZ'}$. 

\section{Extension of the order bound}

We seek to exploit the argument in the order bound a step further by using a partition
\begin{align*}
&C_L(D,G-iP) \backslash C_L(D,G-(i+1)P) \\
=~ &\cup_{j \geq 0} C_L(D,G-iP-jQ) \backslash \\
 &\quad \quad ( C_L(D,G-(i+1)P-jQ) \cup C_L(D,G-iP-(j+1)Q ).
\end{align*}

We apply the argument in the setting of the divisor semigroups.
For a finite set $S$ of rational points, let $\Gamma_S = \cap_{P \in S} \Gamma_P,$ and 
let $\Gamma_S(C) = \{ A \in \Gamma_S : A-C \in \Gamma_S \}$, so that
$\Gamma_S(C) = \cap_{P \in S} \Gamma_P(C).$ Let $\Delta_S(C) = \cup_{P \in S} \Delta_P(C).$

\begin{proposition}
For a divisor class $C$, for a finite set of rational points $S$, and for $P \not \in S,$
\[
\Gamma_S(C) \cap \Gamma_P = \cup_{i \geq 0} \Gamma_{S \cup P}(C+iP).
\]
\end{proposition}
\begin{proof}
Since $P \not \in S,$ $P \in \Gamma_S$, and, for $A-C-iP \in \Gamma_S,$ $A-C \in \Gamma_S.$
Therefore
\[
\Gamma_{S \cup P}(C+iP) = \Gamma_S(C+iP) \cap \Gamma_P(C+iP) \subseteq \Gamma_S(C) \cap \Gamma_P.
\]
For the other inclusion, let $A \in \Gamma_S(C) \cap \Gamma_P.$ Then $A - C \in \Gamma_S$ and $L(A-C) \neq 0$. Choose $i \geq 0 $ 
maximal such that $L(A-C) = L(A-C-iP).$ Then $A-C-iP \in \Gamma_{S \cup P}.$ 
\end{proof}

As a special case $\Gamma(C) \cap \Gamma_P = \cup_{i \geq 0}\, \Gamma_P(C+iP).$ 
Theorem \ref{T:cbdiv} gives a lower bound for $\deg A$, for $A \in \Gamma_P(C).$ Combination of the lower
bounds for $C \in \{ C+iP : i \geq 0 \}$ then gives a lower bound for $\deg A$, for $A \in \Gamma(C) \cap \Gamma_P.$ \\

In combination with $\Gamma(C+iP) \cap \Gamma_Q = \cup_{j \geq 0}\, \Gamma_{\{P,Q\}}(C+iP+jQ),$
we obtain, for $S= \{P,Q\}$, 
\[
\Gamma(C) \cap \Gamma_S = \cup_{i, j \geq 0}\, \Gamma_S(C+iP+jQ).
\]
The next theorem gives a lower bound for $\deg A$, for $A \in \Gamma_S(C).$ For $S=\{P,Q\}$, 
combination of the lower bounds for $C \in \{ C+iP+jQ : i,j \geq 0 \}$ then gives a lower bound for $\deg A$, for 
$A \in \Gamma(C) \cap \Gamma_S.$

\begin{theorem} (Main theorem) \label{T:MTS}
Let $\{ A_1 \leq A_2 \leq \cdots \leq A_w \} \subset \Delta_S(C)$ be
a sequence of divisors with $A_i \in \Delta_{P_i}(C)$, $P_i \in S$, for $i=1,\ldots,w,$
such that $A_{i} - P_{i} \geq A_{i-1}$ for $i=2,\ldots,w$.
Then $\deg A \geq w$, for every divisor $A \in \Gamma_S(C)$ with support disjoint
from $A_w-A_1.$ 
\end{theorem}

\begin{lemma} \label{L:1} 
For $D' \in \Gamma_P(C)$, $\Delta_P(C) \subseteq \Delta_P(D').$
\end{lemma}
\begin{proof}
For $D'-C \in \Gamma_P$, if $A-C \not \in \Gamma_P$ then $A-D' \not \in \Gamma_P$.
\end{proof}

\begin{lemma} \label{L:2}
Let $l_C(A) = l(A) - l(A-C)$. Then
\[
A \in \Delta_P(C) ~\Leftrightarrow~ l_C(A) - l_C(A-P) = 1.  
\]
\end{lemma}
\begin{proof}
\begin{multline*}
l_C(A) - l_C(A-P) = 1 \, 
\Leftrightarrow \, (l(A)-l(A-P)) - (l(A-C)-l(A-C-P)) \,=\, 1 \\
\Leftrightarrow \, l(A)-l(A-P) \,=\, 1 \,\wedge\, l(A-C)-l(A-C-P) \,=\, 0 
\, \Leftrightarrow \, A \in \Delta_P(C).
\end{multline*}
\end{proof}

\begin{proof} \emph{(Theorem \ref{T:MTS})} 
For $A=D' \in \Gamma_S(C) \subseteq \Gamma_{P_i}(C)$ and for 
$A_i \in \Delta_{P_i}(C)$, $A_i \in \Delta_{P_i}(D')$, by Lemma \ref{L:1}, and
$l_{D'}(A_i) = l_{D'}(A_i-P_i)+1$ by Lemma \ref{L:2}. 
With $A_i-P_i \geq A_{i-1}$, there exists a natural map
\[
L(A_{i-1}) / L(A_{i-1}-D') \longrightarrow L(A_i-P_i) / L(A_i-P_i-D').
\]
With $(A_i-P_i)-A_{i-1}$ disjoint form $D'$, the map is injective, since 
$L(A_{i-1}) \cap L(A_i-P_i-D') = L(A_{i-1}-D').$ So that $l_{D'}(A_i-P_i) \geq l_{D'}(A_{i-1})$. 
Iteration over $i$ yields
\[
\deg D' \geq l_{D'}(A_w) \geq l_{D'}(A_{w-1})+1 \geq \cdots \geq l_{D'}(A_1) + w-1  \geq l_{D'}(A_1-P_1) + w.
\]
\end{proof}

To obtain lower bounds with the theorem, we need to construct sequences of divisors in $\Delta_S(C)$.
In the next section we discuss how this can be done effectively.

\section{Efficient computation of the bounds}

The main theorem can be used efficiently to compute coset distances, which in turn can be used to compute two-point code distances. To compute with a certain curve we 
use, for given points $P$ and $Q$, a function $d_{P,Q}$ that encapsulates some geometric properties of the curve \cite{DuuPar08}, see also \cite{BeeNes06JPA}, \cite{Kim94}. 

\begin{lemma}
Let $B$ be a divisor and let $P, Q$ be distinct points. There exists a unique integer $k(B,P,Q)$ such that $B+k'P \in \Gamma_Q$ if and only if $k'\geq k(B,P,Q)$
\end{lemma} 
\begin{proof}
This amounts to showing that $B+kP \in \Gamma_Q$ implies $B+(k+1)P \in \Gamma_Q$. Now use that $\Gamma_Q$ is a semigroup and that $P \in \Gamma_Q.$
\end{proof}
Let us restrict our attention to two-point divisors. Using the previous notation define the following integer valued function.
\begin{align*}
&d_{P,Q}(a)=k(a Q,P,Q)+a
\end{align*}

\begin{theorem}\label{T:d}
For a divisor $A$ with support in $\{P,Q\}$, 
\[
A \in \Gamma_Q ~\Leftrightarrow~ \deg(A)\geq d_{P,Q}(A_Q).
\]
\end{theorem}
\begin{proof}
Let $A=kP+aQ$, then $\deg(A)\geq d_{P,Q}(a)$ is equivalent to $k \geq k(aQ,P,Q)$, which by definition is $aQ+kP \in \Gamma_Q$.
\end{proof}
This property makes the $d$ function a powerful computational tool. Moreover, for $m$ such that $mP\sim mQ$, $d$ is defined modulo $m$. 
In general, the function $d$ depends on the ordering of the points $P$ and $Q$, but it is easy to see that the functions $d_{P,Q}$ and $d_{Q,P}$  satisfy the relation $d_{P,Q}(a) = a+b$ if and only if $d_{Q,P}(b)=a+b$. The function $d = d_{P,Q}$ and the parameter $m$ are enough to 
compute sequences of two-point divisors $A_i \in \Delta_P(C)$ as required by Theorem \ref{T:cbdiv}. A simple application of the Riemann-Roch Theorem give us that we can restrict our search to a finite range of divisors, 
since, for $A \in \Delta_P(C)$, 
\[
\min \{ 0, \deg C \} \;\leq\; \deg A \;\leq\; \max \{ 2g-1, \deg C + 2g-1 \}.
\]
Using $mP \sim mQ$, we can assume moreover, for $A = A_P P + A_Q Q$, that $A_Q \in [0,\ldots,m-1]$. \\

To find long sequences of divisors $A_i \in \Delta_P(C)$ we use a graph theory weight-maximizing algorithm on a rectangular grid $T$, such 
that $T_{i,j}$ is a path of longest length up to the divisor $A$ with degree $i$ and $A_Q=j$ (i.e. $A=iP+j(Q-P)$). \\

\noindent {\bf Computing bounds for the coset $\mathcal{C}(C)\backslash \mathcal{C}(C+P)$ with Theorem \ref{T:cbdiv}.} 
\begin{enumerate}
\item Initialize the first row of $T$ (corresponding to degree $i = \min \{ 0, \deg C \} - 1$) with $0$.
\item Update each row of $T$ successively by the rule $$T_{i,j}=\max\{T_{i-1,j-1},T_{i-1,j}+BP_{i,j}\}$$ where $BP_{i,j}$ is $1$ if $A \in \Delta_P(C)$ and $0$ otherwise,
for the divisor $A$ of degree $i$ with $A_Q=j$.
$BP_{i,j}$ is computed using Theorem \ref{T:d}.
\item Iterate up to the last row (corresponding to degree $i = \max \{ 2g-1, \deg C + 2g-1 \}$).
\item Return the maximum value in the last row.
\end{enumerate}

Using the algorithm we compute bounds for the cosets $\mathcal{C}(C,S)\backslash \mathcal{C}(C+P,S)$ and $\mathcal{C}(C,S)\backslash \mathcal{C}(C+Q,S)$ over all possible divisors $C$. We store them in arrays $CP$ and $CQ$ where the row denotes the degree of $C$ and the column is $C_Q\pmod{m}$. 
For rational points $P$ and $Q$ such that $d_{P,Q} = d_{Q,P}$, we can save some work and obtain the table $CQ$ from $CP$ with the relabeling $CQ_{i,j} = CP_{i,j'}$ for $j' = i-j \pmod{m}.$ 
After computing the coset bounds, 
we traverse all possible coset filtrations of all codes to find bounds for the minimum distances. We use a graph theory flow-maximizing algorithm on a rectangular grid $D$, such that $D_{i,j}$ is a bound for the minimum distance of a code ${\mathcal C}(C)$ with $C$ of degree $i$ and $C_Q=j$.
\\

\noindent {\bf Computing bounds for the distances of all codes ${\mathcal C}(C)$ using $P$-coset and $Q$-coset tables.}

\begin{enumerate}
\item Initialize the last row of $D$ (corresponding to degree $i = 2g$) with $2g$.
\item Update each row of $D$ successively by the rule $$D_{i,j}=\max\{\min\{D_{i+1,j},CP_{i,j}\},\min\{D_{i+1,j+1},CQ_{i,j}\}\}$$
\item Iterate up to the first row (corresponding to degree $i=0$). 
\end{enumerate}

Theorem \ref{T:MTS} can be used to obtain improved lower bounds for the cosets $\mathcal{C}(C)\backslash \mathcal{C}(C+P)$ and
$\mathcal{C}(C)\backslash \mathcal{C}(C+Q).$ The basic algorithm is the same as before with an extra step of using the table of 
all bounds for $\mathcal{C}(C,S) \backslash \cup_{P \in S} \mathcal{C}(C+P,S)$ to produce the $P$- and $Q$-coset tables. \\

\noindent {\bf Computing $P$- and $Q$-coset bounds using Theorem \ref{T:MTS}.} 
\begin{enumerate}
\item Compute table $CS$ with bounds for $\mathcal{C}(C,S)\backslash  \cup_{P \in S} \mathcal{C}(C+P,S)$. This step is done in exactly the same way as the computation for Theorem \ref{T:cbdiv}, but the new table $T$ for each $C$ has the update rule $$T_{i,j}=\max\{T_{i-1,j-1}+BQ_{i,i-j},T_{i-1,j}+BP{i,j}\}$$
\item Initialize $CP$ and $CQ$ at the top row (corresponding to degree $i = 2g$) with $2g$.
\item Compute $CP$ in decreasing row order using the rule $CP_{i,j}=\min\{CP_{i+1,j+1},T_{i,j}\}$.
\item Compute $CQ$ in decreasing row order using the rule $CQ_{i,j}=\min\{CQ_{i+1,j},T_{i,j}\}$.
\end{enumerate}

Once the $P$- and $Q$-coset bounds are known, exactly the same minimizing flow method as the one used for obtaining $d_{DP}$ can be used for an improved bound, denoted $d_{DK}$.

\section{Tables for the Suzuki curve over $\ff_{32}$}




\begin{table}[h] 
\begin{center}
\begin{tabular}{lrrrrrrrrrrrrr}
\toprule
 & & & &\multicolumn{3}{c}{Floor bounds} & &\multicolumn{5}{c}{Order bounds} \\[1ex]
 & &$d_{BPT}$ &~ &$d_{LM}$&$d_{GST}$&$d_{ABZ}$ &~ &$d_{B_0}$&$d_{B}$&$d_{ABZ'}$&$d_{DP}$&$d_{DK}$\\
\midrule
$d_{GOP}$ &  &6352 & &6352 &6352 &6352 & &6352 &6352 &6352 &6352 &6352 \\
$d_{BPT}$ &  &$\cdot$    & &4527 &4551 &4597 & &5260 &5264 &5264 &5264 &5274 \\[1ex]
$d_{LM}$  &  &$\cdot$    & &$\cdot$    &2245 &2852 & &4711 &4729 &4731 &4731 &4757 \\
$d_{GST}$ &  &$\cdot$    & &$\cdot$    &$\cdot$    &2213 & &4711 &4729 &4731 &4731 &4757 \\
$d_{ABZ}$ &  &$\cdot$    & &$\cdot$    &$\cdot$    &$\cdot$    & &4665 &4683 &4685 &4685 &4711 \\[1ex]
$d_{B_0}$ &  &$\cdot$    & &1    &1    &1    & &$\cdot$    &176  &374  &412  &1643 \\
$d_{B}$   &  &$\cdot$    & &1    &1    &1    & &$\cdot$    &$\cdot$    &198  &236  &1565 \\
$d_{ABZ'}$ &  &$\cdot$    & &$\cdot$    &$\cdot$    &$\cdot$    & &$\cdot$    &$\cdot$    &$\cdot$    &38   &1404 \\
$d_{DP}$  &  &$\cdot$    & &$\cdot$    &$\cdot$    &$\cdot$    & &$\cdot$    &$\cdot$    &$\cdot$    &$\cdot$    &1366 \\
\midrule
$d_{GOP}$ &  &1 & &8 &13 &21 & &33 &33 &33 &33 &33 \\
$d_{BPT}$ &  &$\cdot$    & &7 &12 &20 & &32 &32 &32 &32 &32 \\[1ex]
$d_{LM}$  &  &$\cdot$    & &$\cdot$    &7 &15 & &28 &28 &28 &28 &28 \\
$d_{GST}$ &  &$\cdot$    & &$\cdot$    &$\cdot$    &8 & &24 &24 &24 &24 &24 \\
$d_{ABZ}$ &  &$\cdot$    & &$\cdot$    &$\cdot$    &$\cdot$    & &24 &24 &24 &24 &24 \\[1ex]
$d_{B_0}$ &  &$\cdot$    & &1    &1    &1    & &$\cdot$    &1  &5  &5  &6 \\
$d_{B}$   &  &$\cdot$    & &1    &1    &1    & &$\cdot$    &$\cdot$    &5  &5  &6 \\
$d_{ABZ'}$ &  &$\cdot$    & &$\cdot$    &$\cdot$    &$\cdot$    & &$\cdot$    &$\cdot$    &$\cdot$    &1   &6 \\
$d_{DP}$  &  &$\cdot$    & &$\cdot$    &$\cdot$    &$\cdot$    & &$\cdot$    &$\cdot$    &$\cdot$    &$\cdot$    &6 \\
\bottomrule
\end{tabular} 
\end{center}
\caption{Number of improvements of one bound over another (top), and the maximum improvement (bottom), based on $10168$
two-point codes for the Suzuki curve over $\ff_{32}$.}
\label{Table:comparison}
\end{table} 


The Suzuki curve over the field of $q = 2q_0^2$ elements is defined by the equation $y^q+y=x^q_0(x^q+x).$
The curve has $q^2+1$ rational points and genus $g = q_0(q-1)$. 
The semigroup of Weierstrass nongaps at a rational point is
generated by $\{q, q+q_0, q+2q_0, q+2q_0+1\}$. For any two rational points $P$ and $Q$
there exists a function with divisor $(q+2q_0+1)(P-Q).$ Let $m = q+2q_0+1=(q_0+1)^2+{q_0}^2$,
and let $H$ be the divisor class containing $mP \sim mQ$. 
The canonical divisor $K \sim 2(q_0-1)H.$ For the Suzuki curve over the field of $32$ elements
we use $q_0=4, q=32, g=124, m=41, K \sim 6H.$ 
The action of the automorphism group on the rational points of the curve is 2-transitive, so that $d_{P,Q}(a)=d(a)$ does not depend 
on the choice of the points $P$ and $Q$. 
For the Suzuki curve with parameter $q_0$, the $d$ function is given by
\[
d(k) = (q_0-a)(q-1)
\]
where $a$,$b$ are the unique integers such that $|a|+|b| \leq q_0$ and $k \equiv a(q_0+1)+bq_0-q_0(q_0+1) \pmod{m}$. 
A detailed explanation of the geometry behind this result can be found in \cite{DuuPar08}. 
To store the function $d$ as a list, we go through all integers $a$ and $b$ with $|a|+|b| \leq q_0$. To compute $d(k)$ for a single value $k$, we may use
\begin{multline*}
d(k) = (q-1)(2 q_0-q(k-1,2q_0+1)+ \\
-q(r(k-1,2q_0+1),q_0+1)-r(r(k-1,2q_0+1),q_0+1) ),
\end{multline*}
where $q(a,b)$ and $r(a,b)$ are the quotient and the remainder, respectively, when $a$ is divided by $b,$ and $k-1$ is taken modulo $m$. 


Table \ref{Table:comparison} compares the Goppa bound $d_{GOP}$, the base point bound $d_{BPT}$ (an improvement of the Goppa
bound by one whenever $C$ has a base point), the floor bounds $d_{LM}$ \cite[Theorem 3]{LunMcc06}, $d_{GST}$ \cite[Theorem 2.4]{GunStiTas09}, $d_{ABZ}$ \cite[Theorem 2.4]{DuuPar08}, 
and the order bounds $d_{B_0}, d_B, d_{ABZ'}, d_{DP}, d_{DK}.$ 
Each bound was optimized over the full parameter space corresponding to that bound. 
The computations for the Suzuki curve of genus $g=124$ over $\ff_{32}$ were very efficient, computations for the $2g \cdot m = 10168$ 
two-point codes with Goppa designed minimum distance in the range $0, \ldots, 2g-1,$ took less than $10$ minutes on a desktop PC for any given bound. 

As can be seen, the Beelen bound $d_{B}$ offers only slight 
improvement over the weaker Beelen bound $d_{B_0}$. Similar for the improvement
of the Duursma-Park bound $d_{DP}$ over the weaker ABZ bound $d_{ABZ'}$. 
Table \ref{Table:breakdown} gives a further breakdown of the $236$ codes for which $d_{DP}$ improves $d_B$ and the 
$1366$ codes for which $d_{DK}$ improves $d_{DP}$. The improvements are by at most $5$ and $6$, respectively. 

\begin{table}[h]
\begin{center}
\begin{tabular}{lr@{\hspace{2mm}}r@{\hspace{3mm}}r@{\hspace{2mm}}r@{\hspace{2mm}}r@{\hspace{2mm}}r@{\hspace{3mm}}r@{\hspace{4mm}}r@{\hspace{3mm}}r}
\toprule
$d_{DK} - d_{DP}$ &$=$ & &0 &1 &2 &3 &4 &5 &6 \\
\midrule
$d_{DP}-d_{B}$ &$=$ &0 &8603 &656 &356 &198 &50 &6 &63 \\
& &1 &92 &12 &0 &0 &0 &0 &0 \\
& &2 &33 &4 &0 &0 &0 &0 &0 \\
& &3 &74 &4 &1 &0 &0 &0 &0 \\
& &4 &0 &0 &0 &0 &0 &0 &0 \\
& &5 &0 &16 &0 &0 &0 &0 &0 \\
\bottomrule
\end{tabular}
\end{center}
\caption{Improvements of $d_{DP}$ and $d_{DK}$ over $d_B$ for $10168$ two-point codes}
\label{Table:breakdown}
\end{table}

The $63$ codes with $d_{DK}=d_{DP}+6$ all have $d_{DK} = 62$ which agrees with the actual minimum distance (namely realized by
a choice of $62$ points with $P_1 + \ldots + P_{62} \sim 31P + 31 Q$). For each of the $63$ codes,
the coset with $C=23P+23Q$ is the unique coset where Theorem \ref{T:cbdiv} fails to give a lower bound above $56$. For the same coset,
Theorem \ref{T:MTS} gives a lower bound of $62$. \\

Table \ref{Table:comparison} and Table \ref{Table:breakdown} do not show whether the improvements occur for good codes or for poor codes.
For Table \ref{Table:optimal} we select for each degree $\deg C$, i.e. for each given Goppa designed minimum distance, the optimal code with 
respect to each of the bounds $d_{B}, d_{DP}, d_{DK}$ and we compare those. In this case, depending on the degree, the improvements
of $d_{DK}$, obtained with Theorem \ref{T:MTS}, over the bounds $d_B$ and $d_{DP}$ vary between $0$ and $3$.
    


\begin{table}
\begin{tabular}{ccccc}
\begin{tabular}{rr@{\hspace{2mm}}l@{\hspace{2mm}}ll} 
\toprule
$\deg C$ &$d_{DK}$ &$C_{Q}$ & & \\
\midrule
2 &31 &[5] &$\cdot$&$\cdot$ \\
3 &31 & &$\cdot$&$\cdot$ \\
4 &31 & &$\cdot$&$\cdot$ \\
5 &32 &[0] &$\cdot$&$\cdot$ \\
6 &32 & &$\cdot$&$\cdot$ \\
7 &32 & &$\cdot$&$\cdot$ \\
8 &38 &[2] &$\cdot$&3 \\
9 &38 & &$\cdot$&3 \\
10 &38 & &$\cdot$&3 \\
11 &44 &[1] &$\cdot$&$\cdot$ \\
12 &44 &[6] &$\cdot$&$\cdot$ \\
13 &44 & &$\cdot$&$\cdot$ \\
14 &44 & &$\cdot$&$\cdot$ \\
15 &48 &[5] &$\cdot$&1 \\
16 &48 & &$\cdot$&$\cdot$ \\
17 &48 & &$\cdot$&$\cdot$ \\
18 &50 &[8] &$\cdot$&1 \\
19 &50 & &$\cdot$&$\cdot$ \\
20 &51 &[9] &$\cdot$&1 \\
21 &54 &[10] &$\cdot$&3 \\
22 &54 & &$\cdot$&3 \\
23 &54 & &$\cdot$&1 \\
24 &54 & &$\cdot$&3 \\
25 &54 & &$\cdot$&$\cdot$ \\
26 &54 & &$\cdot$&$\cdot$ \\
27 &54 & &$\cdot$&$\cdot$ \\
28 &56 &[12] &$\cdot$&1 \\
29 &56 & &$\cdot$&$\cdot$ \\
30 &58 &[13] &2&2 \\
31 &58 &[10] &2&2 \\
32 &58 & &2&2 \\
33 &59 &[12] &3&3 \\
34 &61 &[13] &3&3 \\
35 &62 &[10 14] &3&3 \\
36 &62 &[17] &3&3 \\
37 &62 & &3&3 \\
38 &62 & &2&2 \\
39 &62 & &2&2 \\
40 &62 & &2&2 \\
41 &62 & &2&2 \\
42 &62 & &2&2 \\
\bottomrule
\end{tabular}
&
\begin{tabular}{rr@{\hspace{2mm}}l@{\hspace{2mm}}ll} 
\toprule
$\deg C$ &$d_{DK}$ &$C_{Q}$ & & \\
\midrule
43 &62 &[5] &$\cdot$&$\cdot$ \\
44 &62 & &$\cdot$&$\cdot$ \\
45 &63 &[5] &$\cdot$&$\cdot$ \\
46 &64 &[0] &$\cdot$&$\cdot$ \\
47 &64 & &$\cdot$&$\cdot$ \\
48 &64 & &$\cdot$&$\cdot$ \\
49 &68 &[2] &$\cdot$&2 \\
50 &68 & &$\cdot$&2 \\
51 &70 &[2, 4] &$\cdot$&2 \\
52 &72 &[1] &$\cdot$&$\cdot$ \\
53 &72 &[4, 6] &$\cdot$&$\cdot$ \\
54 &72 & &$\cdot$&$\cdot$ \\
55 &74 &[7] &2&2  \\
56 &75 &[7] &2&3  \\
57 &75 & &$\cdot$&2 \\
58 &75 & &$\cdot$&$\cdot$ \\
59 &75 & &$\cdot$&$\cdot$ \\
60 &77 &[8] &$\cdot$&$\cdot$ \\
61 &77 & &$\cdot$&$\cdot$ \\
62 &78 &[10] &$\cdot$&$\cdot$ \\
63 &79 &[10] &1&1  \\
64 &80 &[11] &2&2  \\
65 &80 &[8] &$\cdot$&$\cdot$ \\
66 &81 &[5] &1&1  \\
67 &83 &[10] &1&1  \\
68 &85 &[11] &2&2  \\
69 &85 & &1&1  \\
70 &85 &[5, 9] &1&1  \\
71 &87 &[10] &3&3  \\
72 &87 & &2&2  \\
73 &87 &[7, 16] &$\cdot$&$\cdot$ \\
74 &87 & &$\cdot$&$\cdot$ \\
75 &89 &[10] &1&1  \\
76 &91 &[10, 11] &3&3  \\
77 &91 & &1&1  \\
78 &91 &[8] &1&1  \\
79 &91 &[5, 7] &$\cdot$&$\cdot$ \\
80 &93 &[10, 11] &2&2  \\
81 &93 & &2&2  \\
82 &93 & &1&1  \\
83 &93 &[7] &$\cdot$&1 \\
\bottomrule
\end{tabular}
&
\begin{tabular}{rr@{\hspace{2mm}}l@{\hspace{2mm}}ll} 
\toprule
$\deg C$ &$d_{DK}$ &$C_{Q}$ & & \\
\midrule
84 &94 &[1] &$\cdot$&1 \\
85 &94 &[0, 5] &$\cdot$&$\cdot$ \\
86 &96 &[1] &$\cdot$&1 \\
87 &96 &[0] &$\cdot$&$\cdot$ \\
88 &96 & &$\cdot$&$\cdot$ \\
89 &96 & &$\cdot$&$\cdot$ \\
90 &99 &[1, 2] &$\cdot$&$\cdot$ \\
91 &100 &[2, 3] &$\cdot$&$\cdot$ \\
92 &102 &[3] &$\cdot$&$\cdot$ \\
93 &102 & &$\cdot$&$\cdot$ \\
94 &103 &[5] &$\cdot$&$\cdot$ \\
95 &104 &[2] &$\cdot$&$\cdot$ \\
96 &106 &[7] &$\cdot$&$\cdot$ \\
97 &106 & &$\cdot$&$\cdot$ \\
98 &106 & &$\cdot$&$\cdot$ \\
99 &108 &[2, 6] &$\cdot$&$\cdot$ \\
100 &110 &[7] &$\cdot$&$\cdot$ \\
101 &110 & &$\cdot$&$\cdot$ \\
102 &110 & &$\cdot$&$\cdot$ \\
103 &112 &[6] &$\cdot$&$\cdot$ \\
104 &112 &[2] &$\cdot$&$\cdot$ \\
105 &114 &[7] &1&1  \\
106 &114 &[9] &1&1  \\
107 &114 & &$\cdot$&$\cdot$ \\
108 &116 &[2, 6, 7] &$\cdot$&$\cdot$ \\
109 &118 &[7] &2&2  \\
110 &118 & &2&2  \\
111 &118 & &1&1  \\
112 &120 &[6] &1&1  \\
113 &120 &[2] &1&1  \\
114 &121 &[7] &2&2  \\
115 &121 &[2, 6, 9] &1&1  \\
116 &122 &[2, 6, 7, 10, 15] &$\cdot$&$\cdot$ \\
117 &123 &[2, 6, 7] &1&1  \\
118 &124 &[2, 6, 7] &2&2  \\
119 &124 & &1&1  \\
120 &124 &[1, 5, 10, 19] &$\cdot$&$\cdot$ \\
121 &125 &[1] &$\cdot$&$\cdot$ \\
122 &126 &[1] &$\cdot$&$\cdot$ \\
123 &127 &[1] &$\cdot$&$\cdot$ \\
124 &128 &[0] &$\cdot$&$\cdot$ \\
\bottomrule
\end{tabular} 
\end{tabular} \\[1ex]
\caption{Optimal codes (For given $\deg C = 2,\ldots, 124 \,(=g)$, $d_{DK}$ is the maximum lower bound for a two-point code defined with 
$C = C_P P + C_Q Q$, and $C_Q$ gives values for which the maximum is achieved. Suppressed are 
divisors that define subcodes with the same minimum distance as already listed codes. Exchanging $P$ and $Q$ gives a similar
code and listed are only divisors with $C_Q \pmod{m} \leq C_P \pmod{m}.$ The last columns give the amount by which $d_{DK}$ 
exceeds similarly defined maximum lower bounds for $d_{DP}$ and $d_B$, respectively.} 
\label{Table:optimal}
\end{table}


\def\lfhook#1{\setbox0=\hbox{#1}{\ooalign{\hidewidth
  \lower1.5ex\hbox{'}\hidewidth\crcr\unhbox0}}}

\vfill

\noindent Department of Mathematics \\
University of Illinois at Urbana-Champaign\\
1409 W. Green Street (MC-382) \\
Urbana, Illinois 61801-2975\\[1ex]

\noindent \url{http://www.math.uiuc.edu/{~duursma,~rkirov2}}

\end{document}